\documentclass[12 pt]{article}
\usepackage{amscd, amssymb,amsmath,amsthm}
\newtheorem{Relation}{Relation}
\newtheorem{Lemma}{Lemma}
\newtheorem{Theorem}{Theorem}
\newtheorem{Corollary}{Corollary}

\newcommand{\ev}{\operatorname{ev}}

\newcommand{\proj}{\mathbb{P}}
\newcommand{\Z}{\mathbb{Z}}

\newcommand{\rarr}{\rightarrow}
\newcommand{\oh}{{\mathcal{O}}}
\newcommand{\com}{\mathbb{C}}
\newcommand{\Q}{\mathbb{Q}}

\newcommand{\Lan}{\Big \langle}
\newcommand{\Ran}{\Big \rangle}
\newcommand{\lan}{\langle}
\newcommand{\ran}{\rangle}

\usepackage{leqno}
\begin{document}

\title{A topological view of Gromov-Witten theory}
\author{D.~Maulik and R.~Pandharipande}
\date{March 2005}
\maketitle

\setcounter{section}{-1}
\section{Introduction}
\subsection{Overview}
Let $V$ be a nonsingular,  
 complex, projective variety containing a nonsingular divisor $W$.
The {\em absolute} Gromov-Witten theory of $V$ is defined by
 integrating descendent classes over the moduli space of stable
maps to $V$.
The {\em relative} Gromov-Witten theory of the pair $(V,W)$ is defined by descendent integration 
over the
space of stable relative maps to $V$ with prescribed tangency data along $W$.

We present here a systematic study of relative Gromov-Witten theory 
via universal relations.
We find the relative theory does {\em not} provide
new invariants: the relative theory is completely determined by the absolute
theory. The relation between the relative and absolute theories is guided by a 
strong analogy to classical topology.

Our results open new directions in the subject.
For example, we present a complete mathematical determination of the
Gromov-Witten theory (in all genera) of the Calabi-Yau quintic hypersurface in $\proj^4$.

\subsection{Leray-Hirsch}
\label{ov}
Let $X$ be a nonsingular, complex, projective variety 
equipped with a line bundle $L$.
Let $Y$ be the projective bundle
$\proj(L \oplus \oh_X)$, and
let $\pi$ be the projection map,
$$\pi:Y \rarr X.$$
The summands $L$ and $\oh_X$ respectively determine divisors
$$D_0, D_\infty \subset Y$$
isomorphic to $X$ via $\pi$.

We consider four descendent Gromov-Witten 
theories (in all genera) of the projective bundle 
$Y$:
the absolute Gromov-Witten theory of $Y$ and the
relative Gromov-Witten theories of the three pairs
$$(Y,D_0), \ (Y,D_\infty), \ (Y, D_0 \cup D_\infty).$$
The first result of the paper is a reconstruction theorem for the Gromov-Witten theories
of $Y$ in terms of $X$.

\begin{Theorem}\label{mn} All four  theories  of $Y$ can be uniquely and effectively 
reconstructed from
the Gromov-Witten theory of $X$ and the class $c_1(L)\in H^2(X,\Q)$.
\end{Theorem}

Theorem \ref{mn} is proven in Section \ref{aaa} by exhibiting an explicit 
set of recursions. The two main
techniques used are localization (with respect to the natural fiberwise $\com^*$-action 
 on $Y$) and degeneration. Equivariant relative
invariants of $\proj^1$ appear as constants in the recursions. 
We view Theorem \ref{mn} as a Leray-Hirsch
result in Gromov-Witten theory.

\subsection{Relative in terms of absolute}
Let $V$ be a nonsingular,  
 complex, projective variety containing a nonsingular divisor $W$.
Let 
 $$N\rarr W$$
be the normal bundle of $W$ in $V$.

Let $\mathcal{F}$ be the degeneration to the normal cone of $W$,
 the
 blow-up of $V\times \com$ along the subvariety
$W\times 0$.
Let $$\epsilon: \mathcal{F} \rarr \com$$
be the projection to the second factor.
We find
$$\epsilon^{-1}(0) = V \cup_{W} \proj(N \oplus \oh_W)$$
where the inclusion
$$W\subset \proj(N \oplus \oh_W)$$
is determined by summand $N$.
The degeneration formula \cite{EGH,IP,LR,L} expresses the absolute Gromov-Witten theory of
$V$ in terms of the relative theories of the pairs
$(V,W)$ and $(\proj(N\oplus \oh_W), W)$.

Theorem \ref{mn} together with an inversion of the degeneration formula yields
the following result proven in Section \ref{bbb}.

\begin{Theorem}
\label{wd}
The relative Gromov-Witten theory of the pair $(V,W)$ can be uniquely and effectively
reconstructed from the 
absolute theory of $V$, the absolute theory of $W$, and the restriction map
$H^*(V,\Q) \rarr H^*(W,\Q)$.
\end{Theorem}

We view Theorem \ref{wd} as a Gromov-Witten analogue of the standard long exact 
sequence relating
absolute and relative cohomology theories.

\subsection{Mayer-Vietoris}
Let $V$ be a nonsingular, complex, projective variety. Let
$$\epsilon:{\mathcal V} \rarr \Delta$$
be a flat family over a disk $\Delta \subset \com$ at the origin
satisfying:
\begin{enumerate}
\item[(i)] $\mathcal {V}$ is nonsingular,
\item[(ii)] $\epsilon$ is smooth over the punctured disk $\Delta^*=\Delta\setminus \{0\}$,
\item[(iii)] $\epsilon^{-1}(1)\stackrel{\sim}{=} V$,
\item[(iv)] $\epsilon^{-1}(0)=V_1\cup_W V_2$ is a normal crossings divisor in ${\mathcal V}$.
\end{enumerate}
The family $\epsilon$ defines a canonical map
$$H^*(V_1\cup_W V_2,\Q) \rarr H^*(V,\Q)$$
with image defined to be the {\em nonvanishing} cohomology of $V$.

The degeneration formula and Theorem \ref{wd} together
yield a Mayer-Vietoris result.

\begin{Theorem}\label{mv}
The Gromov-Witten theory of the nonvanishing cohomology of
$V$ can be uniquely and effectively
reconstructed from the absolute theories of $V_1$, $V_2$, and $W$ and the
restriction maps
$$H^*(V_1,\Q) \rarr H^*(W,\Q), \ \ H^*(V_2,\Q)\rarr H^*(W, \Q).$$ 
\end{Theorem}

\subsection{Hypersurfaces} \label{sch}
\subsubsection{Hypersurface pairs}

A {\em hypersurface pair} $(V,W)$ is a nonsingular hypersurface
$V\subset \proj^r$ together with a nonsingular divisor $W\subset V$ 
defined by a complete intersection in $\proj^r$.

Let $\beta\in H_2(V,{\mathbb Z})$ be a curve class.
Let $\stackrel{\rightarrow}{\mu}$ be an ordered partition,
$$\sum_j \mu_j = \int_\beta[W],$$
with positive parts.
The moduli space $\overline{M}_{g,n}(V/W,\beta, \stackrel{\rightarrow}{\mu})$ parameterizes
stable relative maps from genus $g$, $n$-pointed curves to $V$ of class $\beta$
with 
 multiplicities along $W$ determined by $\stackrel{\rightarrow}{\mu}$.

The relative conditions in the theory correspond to partitions {\em weighted} by
the cohomology of $W$.
Let  $\delta_1, \ldots, \delta_{m_W}$ be a basis of $H^*(W,{\mathbb Q})$.
A cohomology weighted partition ${\nu}$
consists of an {\em unordered} set of pairs,
$$\left\{ (\nu_1, \delta_{s_1}), \ldots, (\nu_{\ell(\nu)}, \delta_{s_{\ell(\nu)}}) 
\right\}, $$
where $\sum_j \nu_j$ is an {\em unordered} partition of $\int_\beta [W]$.
The automorphism group, $\text{Aut}(\mathbf{\nu})$, consists of 
permutation symmetries of ${\mathbf \nu}$.

The {\em standard} order on the parts of $\nu$ is
$$(\nu_i,\delta_{s_i})> (\nu_{i'},\delta_{s_{i'}})$$
if $\nu_i>\nu_{i'}$ or if $\nu_i=\nu_{i'}$ and $s_i>s_{i'}$.
Let $\stackrel{\rarr}{\nu}$ denote the partition $(\nu_1,\ldots, \nu_{\ell(\nu)})$
obtained from the standard order.

The descendent Gromov-Witten invariants of 
the hypersurface pair are defined by integration against the
virtual class of the moduli of maps. 
Let $\gamma_1,\ldots, \gamma_{m_V}$ be a basis of $H^*(V,\Q)$, and
let
\begin{multline*}
\left. \Lan  \tau_{k_1}(\gamma_{l_1}) \cdots
\tau_{k_n}(\gamma_{l_n})\ \right| {\mathbf \nu} 
\Ran^{V/W}_{g,\beta} = \\
 \frac{1}{|\text{Aut}(\nu)|}
\int_{[\overline{M}_{g,n}(V/W,\beta,\stackrel{\rarr}{\nu})]^{vir}} 
\prod_{i=1}^n \psi_i^{k_i}\text{ev}_i^*(\gamma_{l_i}) \cup 
\prod_{j=1}^{\ell(\nu)} \text{ev}^*_j(\delta_{s_j}).
\end{multline*}
Here, the second evaluations, 
$$\text{ev}_j: \overline{M}_{g,n}(V/W,\beta,\stackrel{\rarr}\nu) \rarr W.$$
are determined by the relative points.

Gromov-Witten invariants are defined (up to sign) for {\em unordered}
weighted partitions ${\nu}$. To fix the sign, the integrand
on the right side requires an ordering. The ordering is corrected
by the automorphism prefactor.

\subsubsection{Simple classes}

A class $\gamma\in H^*(V,\Q)$ is {\em simple} if $\gamma$ lies in the image of
the restriction map
$$H^*(\proj^r,\Q) \rarr H^*(V,\Q).$$
The  simple Gromov-Witten theory of $V$ consists of 
the integrals of descendents of simple classes. Similarly, the
simple Gromov-Witten theory of the pair $(V,W)$ consists of integrals of
descendents of simple classes with {\em no restrictions on the 
cohomology classes of $W$ in the relative constraints}.

A refinement of Theorem  \ref{wd} proven in Section \ref{bbb} 
is valid for the geometry of the hypersurface pair
$(V,W)$.

\begin{Corollary}\label{rwd}
The simple Gromov-Witten theory of a hypersurface pair $(V,W)$ can be uniquely and effectively
reconstructed from  
the simple 
Gromov-Witten theory of $V$,  the full Gromov-Witten theory of $W$, and the restriction map
$H^*(V,\Q) \rarr H^*(W,\Q).$  
\end{Corollary}

\subsubsection{Curves, surfaces, and 3-folds}
Nonsingular curves have a rich Gromov-Witten theory 
including descendents of odd classes. The Gromov-Witten theory of
curves is fully determined in \cite{OP1,OP2,OP3}.

Nonsingular surfaces in $\proj^3$ of degree up to 3 are rational with Gromov-Witten theories
determined by localization \cite{giv,GP}. The $K3$ surface in degree 4 is
holomorphic symplectic --- hence all
Gromov-Witten invariants of the $K3$ vanish for nonconstant maps \cite{bl}. 
We will present a complete
scheme for calculating the simple Gromov-Witten theory of surfaces of degree 5 and higher.

Nonsingular 3-folds in $\proj^4$ are determined by 
localization methods only in degrees 1 and 2. We present a complete scheme for 
calculating the simple Gromov-Witten theory of hypersurfaces
of degree 3, 4, and 5 in $\proj^4$. 

Gromov-Witten theory is not interesting for
3-fold hypersurfaces of degree greater than 5
since the moduli spaces of nonconstant maps have negative dimension.

\subsubsection{Calculation scheme}
Let $X\subset \proj^r$ be a generic hypersurface of degree $d \geq 2$ with equation $f$.
Let 
$$h=h_1h_2$$
be a product of generic polynomials of degree $d_1-1$ and $1$.
Let $t$ be a coordinate on $\com$.
The ideal
$$(tf-h)$$
determines a subvariety 
$${\mathcal X} \subset \proj^r \times \com$$
flat over 
over $\com$.
The generic element of the family is a nonsingular hypersurface. 
The fiber over $0\in \com$ is a union,
$${\mathcal X}_0= X_1 \cup_I X_2,$$
of a hypersurface $X_1$ with equation $h_1$
and a hyperplane $X_2$ with equation $h_2$ 
along a complete
intersection $I\subset \proj^r$ of type $(d_1-1,1)$.

The total space of ${\mathcal X}$ is singular. The singular locus of ${\mathcal X}$ over
$0\in \com$ is a complete intersection $S\subset \proj^r$ of type $(d_1,d_1-1,1)$
contained in $I$. Locally analytically, the singularities of ${\mathcal X}$ over $0\in \com$ are
translates of the 3-fold double point.

Let $\widetilde{\mathcal X}$ denote the blow-up of ${\mathcal X}$ along the
Weil divisor $X_2$. The total space $\widetilde{\mathcal X}$ is nonsingular
over $0\in \com$. 
The fiber over $0\in \com$ is a union,
$$\widetilde{\mathcal X}_0= X_1 \cup_I \widetilde{X}_2,$$
where $X_1$ and $I$ are as before and
$\widetilde{X}_2$ is the blow-up of $X_2$ along $S$. 

By the degeneration formula, the simple Gromov-Witten theory of $X$ is determined
by the simple Gromov-Witten theories of the pairs 
$(X_1,I)$ and $(\widetilde{X}_2,I)$. By Corollary \ref{rwd},
the simple Gromov-Witten theory of $(X_1,I)$ is determined by 
the simple theory of
$X_1$ and the full theory of $I$. The simple Gromov-Witten theory of
$(\widetilde{X}_2,I)$ requires, in  addition, the simple theory of $\widetilde{X}_2$.

By the application of Lemma \ref{ss} of Section \ref{ccc} below to $\widetilde{X}_2$,
 we conclude
the simple Gromov-Witten theory of $X$ is determined by the simple theories
of the hypersurfaces
$$X_1, X_2\subset \proj^r$$
of lower degree and the full theories of
the varieties $I$ and $S$ of lower dimension.

\subsubsection{Blow-up Lemma}
\label{ccc}
Let $V$ be a nonsingular, projective variety. Let
$Z\subset V$ be the nonsingular complete intersection of two nonsingular divisors
$$W_1,W_2\subset V,$$
and let $\widetilde{V}$ be the blow-up of $V$ along $Z$.

\begin{Lemma}\label{ss}
The Gromov-Witten theory of $\widetilde{V}$ is
uniquely and effectively determined  by
the Gromov-Witten theories of $V$, $W_1$, and $Z$
and the restriction maps
$$H^*(V,\Q) \rarr H^*(W_1,\Q) \rarr H^*(Z,\Q).$$ 
\end{Lemma}

Only the absolute Gromov-Witten theory of {\em one} of the divisors
is needed in Lemma \ref{ss}. However, an optimal result should
avoid both divisors. Lemma \ref{ss} is proven in Section \ref{edd}.

\subsubsection{Surfaces}
The calculation scheme determines the simple Gromov-Witten invariants of 
surfaces of degree $d\geq 5$ in $\proj^3$.

For example, the simple Gromov-Witten theory of the degree 5 surface $S_5$
is determined in 
terms of the Gromov-Witten theories of the following spaces:
 $$\proj^2, \ S_4, \ C_4, \ \proj^0,$$
where $S_4$ is a $K3$ surface and
$C_4$ is a nonsingular quartic plane curve.

For surfaces of general type, Gromov-Witten invariants in the adjunction genus
with
{\em primary} field insertions are determined by gauge theory: Taubes' tetrology 
connects these Gromov-Witten invariants to Seiberg-Witten theory \cite{t1,t2,t3,t4}.
In Section \ref{exx}, we present a calculation of
$$\langle 1\rangle^{S_5}_{6,K} = -1$$
by our method.
The structure of Gromov-Witten theory in other genera or
with descendent insertions is not known.

\subsubsection{3-folds}
The calculation scheme determines the simple Gromov-Witten invariants of 
hypersurfaces of degree $3$, $4$, and $5$ in $\proj^4$ in terms of
known theories. 


The Calabi-Yau
quintic 3-fold $Q\subset \proj^4$ appears to be the most difficult hypersurface captured
by the scheme. The Gromov-Witten invariants of $Q$ are 
determined in terms of the Gromov-Witten theories of the following spaces:
$$\proj^3,\ \proj^2,  \ S_2,\ S_3,\ S_4, \ C_{1,2}, \ C_{2,3}, \ C_{3,4}, \ C_{4,5}.$$
Here, $S_d\subset \proj^3$ is a nonsingular degree $d$ surface, and
$C_{d_1,d_2}\subset \proj^3$ is a nonsingular complete intersection curve of type $(d_1,d_2)$.

The quintic scheme, the first mathematical determination of the
Gromov-Witten theory of $Q$, is
presented in Section \ref{zzwe} 

\subsubsection{Further directions}
Our calculation scheme using degeneration and Mayer-Vietoris can be
pursued in several other contexts --- hypersurfaces are treated here
as the first illustration of the method. For example, the simple
Gromov-Witten theories of all
surface and 3-fold complete intersections in projective space are determined
similarly.

\subsection{Gathmann's proposal}

Our last topic concerns Gathmann's proposal for the calculation of higher genus 
Gromov-Witten 
invariants of the quintic $Q\subset \proj^4$.

Gathmann has studied the Gromov-Witten invariants of $Q$ in genus 0 and 1 via a 
relation to
the Gromov-Witten theory of $\proj^4$ --- the opposite direction of the
scheme discussed in Section \ref{sch}. Theorem \ref{mn} allows us to pursue 
Gathmann's proposal in all genera.

Let $\mathcal{G}$ be the blow-up of $\proj^4\times \com$ along the subvariety
$Q\times 0$.
Let $$\epsilon: \mathcal{G} \rarr \com$$
be the projection to the second factor.
We find
$$\epsilon^{-1}(0) = \proj^4 \cup_{Q} \proj(\oh_Q(5) \oplus \oh_Q)$$
where the inclusion
$$Q\subset \proj(\oh_Q(5)\oplus \oh_Q)$$
is determined by summand $\oh_Q(5)$.

The degeneration formula \cite{EGH,IP,LR,L} expresses the absolute Gromov-Witten 
theory of
$\proj^4$ in terms of the relative theories of the pairs
$(\proj^4,Q)$ and $(\proj(\oh_Q(5)\oplus \oh_Q), D_0)$.
The absolute Gromov-Witten theory of $\proj^4$ may be computed by several
methods \cite{g3,giv,GP}. Theorem \ref{mn} reduces the relative theory of
$(\proj(\oh_Q(5)\oplus \oh_Q), D_0)$ to the absolute theory of $Q$.
The degeneration formula then provides a system of equations for
the relative invariants of $(\proj^4,Q)$ {\em and} the Gromov-Witten invariants
$N_{g,d}$ of $Q$.

\vspace{+10pt}
\noindent{\bf Conjecture 1.} {\em The system of equations obtained from the
degeneration formula and Theorem \ref{mn} determines both the
relative theory of the pair $(\proj^4, Q)$  and the Gromov-Witten invariants
$N_{g,d}$ of $Q$.}
\vspace{+10pt}

Gathmann pursued the above method in genus 0 and 1 \cite{g1,g2,g4} via a different
 approach to the
reduction of the relative theory of the pair $$(\proj(\oh_Q(5)\oplus \oh_Q), D_0)$$ to the
absolute theory of $Q$.
Theorem \ref{mn}, however, is valid for {\em all} genera. Gathmann's method and 
Theorem \ref{mn}
together provide a computation scheme for all $N_{g,d}$ so long as the
arising equations are nonsingular. 

We have proven Conjecture 1 in genus 2. However,
even if the nonsingularity is not  known
beforehand, the computation can be undertaken. 
Gathmann's proposal,
though difficult and not yet certain for $g\geq 3$, appears more suitable for 
calculations
than the complete quintic scheme discussed in Section \ref{sch}.

A very different approach to the Gromov-Witten theory of $Q$ in genus 1 has been
advanced in a series of papers by Zinger and collaborators \cite{lz,vz,z1,z2}.

\subsection{Acknowledgments}
We thank J. Bryan for many helpful discussions. In particular, the
quintic surface calculation in Section \ref{exx} was motivated by conversations with him. 
We also benefitted from discussions with
T. Graber, E. Katz, A. Okounkov,  
and A. Zinger. 

D.~M. was partially
supported by an NSF graduate fellowship. R.~P. was partially supported by the NSF and
the Packard foundation.

\section{Leray-Hirsch} \label{aaa}
\subsection{Notation}
\subsubsection{Cohomology}
Let $X$ be a nonsingular, projective variety 
equipped with a line bundle $L$, and let $Y$ be the
projective bundle
$$\pi:\proj(L \oplus \oh_X)\rarr X$$
with sections $D_0,$ $D_\infty$ corresponding to
the summands $L$, $\oh_X$ respectively.

Let $\delta_{1}, \dots, \delta_{m_X}$ be a basis of $H^*(X,\Q)$ containing 
the identity element. We will often denote the identity by $\delta_{Id}$.
The {\em degree} of $\delta_i$ is the real grading in $H^*(X,\Q)$.
We  view $\delta_i$ as an element of $H^*(Y,\Q)$ via pull-back by $\pi$.

Let 
$[D_{0}],
[D_{\infty}] \in H^{2}(Y,\Q)$
denote the cohomology classes associated to the divisors.
Define classes in  $H^*(Y,\Q)$ by
\begin{eqnarray*}
\gamma_{i} & = &  \delta_{i}\ , \\
\gamma_{m_X+i} & = & \delta_i \cdot [D_0]\ , \\
\gamma_{2m_X+i} & = & \delta_i \cdot [D_\infty]\ .
\end{eqnarray*}
We will use the following notation:
\begin{eqnarray*}
 \gamma^{\delta}_i & = &  \delta_{i \ \text{mod}\ m_X}\ , \\
\gamma^{D}_i & = & 1, [D_0],\ \text{or} \  [D_\infty]\ .
\end{eqnarray*}
The second assignment depends upon the integer part of $i/m_X$.
The set $\{\gamma_1,\ldots, \gamma_{2m_X}\}$ determines a basis of $H^*(Y,\Q)$.

\subsubsection{Theorem \ref{mn} for the absolute theory of $Y$} \label{fca}
There is a fiberwise $\com^*$-action on $Y$ determined by scaling
the second factor in the sum $L \oplus \oh_X$.
The absolute theory of $Y$ can be directly computed via
the virtual localization formula of \cite{GP}.

The $\com^*$-action on $Y$ induces a canonical $\com^*$-action on the
moduli space of stable maps $\overline{M}_{g,n}(Y,\beta)$.
The $\com^*$-fixed loci of $\overline{M}_{g,n}(Y,\beta)$  
are determined by bipartite graphs. The vertices correspond to
spaces of stable maps to $D_0$ or $D_\infty$ --- both targets are
canonically isomorphic to $X$.
The virtual localization formula reduces the Gromov-Witten invariants
of $Y$ to Hodge integrals in the Gromov-Witten theory of $X$.
The Hodge insertions may be removed by the relations of \cite{FP}.
The proof of Theorem \ref{mn} for the absolute theory of $Y$ is complete.

\subsubsection{Brackets} 
We will use the following bracket notation for the Gromov-Witten invariants of
the pair $(Y,D_0)$:
\begin{multline*}
\Lan {\mu}\Big| \prod_i \tau_{k_{i}}(\gamma_{l_{i}})\Ran_{g,\beta} = \\
\frac{1}{\mathrm{Aut}(\mu)}\int_{[\overline{M}_{g,n}(Y/D_{0},
\beta,\stackrel{\rarr}{\mu})]^{\mathrm{vir}}}
\prod_i \psi_{i}^{k_{i}}\ev^{\ast}(\gamma_{l_{i}})
\cup\prod_j \ev^{\ast}(\delta_{r_{j}}),
\end{multline*}
where 
$$\mu=\{ (\mu_1, \delta_{r_1}), \ldots, (\mu_{\ell(\mu)}, \delta_{r_{\ell(\mu)}})\}$$
is a  partition weighted by the cohomology of $X$ and 
$$\sum_j \mu_j=\int_\beta [D_0].$$
Relative invariants are defined only when $\int_\beta[D_0]\geq 0$.

For the pair $(Y,D_\infty)$, the relative conditions will be written on the 
right side of
the bracket:
$$ \Lan \prod_i \tau_{k_{i}}(\gamma_{l_{i}})\Big|{\nu}\Ran _{g,\beta}.$$  
The invariants of the pairs $(Y, D_0)$ and $(Y,D_\infty)$
are  termed {\em type I}.

For the pair $(Y,D_0\cup D_\infty)$, the relative conditions for $D_0$ 
will be written on the  left side and the relative conditions for $D_\infty$ will be 
written on the right side:
$$ \Lan \mu \Big| \prod_i \tau_{k_{i}}(\gamma_{l_{i}})\Big|{\nu}\Ran_{g,\beta}.$$  
The invariants of $(Y,D_0\cup D_\infty)$
are termed {\em type II}.

The above brackets denote Gromov-Witten invariants with {\em connected}
domain curves. Disconnected invariants arise naturally in 
the degeneration formula. 
We will treat disconnected invariants as products of connected invariants
except in the study of rubber targets in Section \ref{rc}. However, our proof
of Theorem \ref{mn} is valid without assuming the product rule.
The connected/disconnected issue will be discussed carefully in Section \ref{disco}.

\subsubsection{Partition terminology}
The following constants associated to a weighted partition $\mu$ will arise often:
\begin{enumerate}
\item[$\bullet$] $\deg(\mu)= \sum_i \deg(\delta_{r_i})$, the total degree of the
cohomology weights,
\item[$\bullet$] $Id(\mu)$ equals the number occurences of the pair $(1,\delta_{Id})$ in $\mu$,
\item[$\bullet$]
$\mathfrak{z}(\mu)= \prod_i \mu_i \cdot | \text{Aut}(\mu)|.$
\end{enumerate}

We assume the cohomology basis $\delta_1, \ldots, \delta_{m_X}$ is self dual
with respect to the Poincar\'e pairing. Then, to each wieghted partition $\mu$, a
dual partition $\mu^\vee$ is defined by taking the Poincar\'e duals of the
cohomology weights.

\subsubsection{Orderings}
All Gromov-Witten invariants $\lan,\ran_{g,\beta}$ vanish if $$\beta\in H_2(Y,\Z)$$ is not
an effective curve class. We define a partial ordering on $H_2(Y,\Z)$
as follows:
$$\beta'< \beta$$
if $\beta-\beta'$ is a nonzero effective curve class.

The set of pairs $(m,\delta)$ where $m\in \Z_{>0}$ and
$\delta\in H^*(X,\Q)$ is partially ordered by the following
{\em size} relation
\begin{equation}\label{ssf}
(m,\delta) >  (m',\delta')
\end{equation}
if $m> m'$ or if $m=m'$ and $\deg(\delta)>\deg(\delta')$.

Let $\mu$ be a partition weighted by the cohomology of $X$,
$$\mu=\{ (\mu_1, \delta_{r_1}), \ldots, (\mu_{\ell(\mu)}, \delta_{r_{\ell(\mu)}})\}.$$
We may place the pairs of $\mu$ in decreasing order by
size \eqref{ssf}.
A {\em lexicographic} ordering on weighted partitions
is defined  as follows:
$$\mu \stackrel{l}{>} \mu'$$
if, after placing $\mu$ and $\mu'$ in decreasing order by size,
the first pair
 for which $\mu$ and $\mu'$ differ in size is larger for $\mu$.

\subsection{Fiber classes}\label{fclass}
Let $[F]\in H_2(Y,\Z)$ denote the class of a fiber of $\pi$. The {\em fiber class} 
invariants of
type I and II 
are those for which $\beta$ is a (possibly zero) multiple of $[F]$.
Our first goal is to calculate the fiber class invariants of both types in terms of the
classical cohomology of $X$.

Consider a connected type I invariant of a fiber class,
\begin{equation}\label{exex}
\Lan \mu \Big| \prod_{i} \tau_{k_i}(\gamma_{l_i}) \Ran_{g,d[F]}.
\end{equation}
We determine the fiber class invariant \eqref{exex}  
from the equivariant theory of $\proj^{1}$.
 
The moduli space of stable relative maps 
 $$\overline{M}_{Y}=\overline{M}_{g,n}(Y/D_{0},d[F],\stackrel{\rarr}{\mu}),$$ 
is fibered over $X$,
$$\pi: \overline{M}_Y \rarr X,$$
with fiber isomorphic to the moduli space of maps to $\proj^1$ relative to 0,
$$\overline{M}_{\proj^{1}}=\overline{M}_{g,n}(\proj^{1}/0,d,\stackrel{\rarr}{\mu}).$$
 In fact, $\overline{M}_{Y}$ is the fiber bundle constructed from the
principal $\com^{*}$-bundle 
associated to $L$ and a standard $\com^{*}$-action on $\overline{M}_{\proj^{1}}$.  

The $\pi$-relative
obstruction theory of $\overline{M}_{Y}$ is obtained 
from the $\overline{M}_{\proj^{1}}$-fiber bundle structure over $X$.  
The relationship between the $\pi$-relative virtual fundamental class
$[\overline{M}_{Y}]^{vir_\pi}$ and the 
virtual fundamental class $[\overline{M}_Y]^{vir}$ is given by the equation
\begin{equation}\label{fgv}
[\overline{M}_{Y}]^{{vir}} 
=  c_{{top}}(\mathbb{E}\boxtimes T_{X}) \cap [\overline{M}_{Y}]^{vir_\pi} 
\end{equation}
where $\mathbb E$ is the Hodge bundle.
We may rewrite \eqref{fgv} as
\begin{equation*}
[\overline{M}_{Y}]^{{vir}} 
=  \sum_q h_q\Big(c_1(\mathbb E), c_2(\mathbb E), \ldots\Big)\ t_q\Big(c_1(T_X), c_2(T_X), 
\ldots\Big)
\ \cap [\overline{M}_{Y}]^{vir_\pi} 
\end{equation*}
where $h_q$ and $t_q$ are polynomials.

The invariant \eqref{exex} can then be computed by pairing  
cohomology classes in $X$ with the results of equivariant integrations in
the Gromov-Witten theory of $\overline{M}_{\proj^1}$.
We write 
\begin{multline}
\Lan \mu \Big| \prod_{i} \tau_{k_i}(\gamma_{l_i}) \Ran_{g,d[F]} =\\
\frac{1}{\mathrm{Aut}(\mu)} \sum_q 
\int_X \Big( t_q\prod_i \gamma^\delta_{l_i} \prod_j \delta_{r_j}  \cap
\pi_*\big(h_q\prod_i \psi_i^{k_i} \ev_i^*(\gamma_{l_i}^D) \cap [\overline{M}_Y]^{vir_\pi}\big) 
\Big)
\end{multline}
The interior push-forward 
$$\pi_*\big(h_q\prod_i \psi_i^{k_i} 
\ev_i^*(\gamma_{l_i}^D) \cap [\overline{M}_Y]^{vir_\pi}\big)$$
is obtained from
the corresponding  Hodge integral in the equivariant Gromov-Witten theory
of $\proj^1/0$  after
replacing the hyperplane class on $\com\proj^\infty$  by $c_{1}(L)$.

The argument for the pairs $(Y,D_\infty)$ and $(Y,D_0\cup D_\infty)$ is
identical. 
The required Hodge integrals in the 
equivariant relative Gromov-Witten theory of $\proj^1$ are
fully determined by the methods of the papers \cite{FP,FP2,OP1,OP2}.

\subsection{Distinguished type II invariants} \label{xxw}

{\em Distinguished} invariants of type II are integrals for which there
 is a distinguished marked point $p$ with a {\em pure} 
 cohomology condition of the form 
$[D_0]\cdot \delta$ where
$$\deg(\delta) >0.$$
We will write distinguished invariants as
$$\Lan \mu \Big| \tau_0([D_0] \cdot \delta)\cdot \omega \Big|\nu \Ran_{g,\beta}$$
where $\tau_0([D_0]\cdot \delta)$ is the distinguished insertion and $\omega$
denotes the product of the non-distinguished insertions.
Let $\|\omega\|$ denote the number of  
non-distinguished insertions.

We will compute distinguished type II invariants by an inductive algorithm. 
A partial ordering $\stackrel\circ{<}$
on the set of distinguished type II invariants is defined as follows: 
$$\Lan \mu'\Big|
\tau_0([D_0]\cdot \delta')\cdot \omega'\Big |\nu'\Ran_{g',\beta'} \stackrel{\circ}{<} 
\Lan \mu \Big|\tau_{0}([D_0]\cdot \delta)\cdot  \omega
\Big|\nu \Ran_{g,\beta}$$
if one of the conditions below holds
\begin{enumerate}
\item[(1)] $\beta' < \beta$, 
\item[(2)] equality in (1) and $g' < g$,
\item[(3)] equality in (1-2) and $\| \omega'\| < \| \omega \|$,
\item[(4)] equality in (1-3) and $\deg({\mu}') > \deg({\mu})$,
\item[(5)] equality in (1-4) and $\deg({\nu}') > \deg({\nu})$,
\item[(6)] equality in (1-5) and $\deg(\delta')> \deg(\delta)$,
\item[(7)] equality in (1-6) and ${\nu}' \stackrel{l}{>} {\nu}$,
\end{enumerate}

For any given distinguished invariant of type II , there are only finitely many 
distinguished invariants of type II lower in the partial ordering.  
Our algorithm consists of relations between 
type I and type II invariants which allow us to move down 
the partial ordering.

\subsection{Relation \ref{A}}\label{ax1}
Relation \ref{A} expresses a distinguished invariant  of type II in terms of 
type I invariants and strictly lower distinguished invariants of  type II with
respect to $\stackrel{\circ}{<}$.

Fix $g$ and $\beta>0$. Type I and II invariants of genus $g$ and class $\beta$
will be viewed as {\em principal} terms of the equations below.
All type I and II invariants 
of $Y$ with $$\beta^{\prime} < \beta$$ or 
$$\beta^{\prime} = \beta\ {\text {and}}\ g^{\prime}<g$$
are viewed as {\em non-principal} terms.
The non-principal terms are inductively determined and, therefore, omitted
in the equations.

Let $R$ denote the distinguished type II invariant 
$$\Lan {\mu}\Big|\tau_0([D_0]\cdot \delta) \cdot \omega\Big| {\nu}\Ran_{g,\beta}$$
with $\deg(\delta)>0$
and relative conditions
$$\mu=\big\{ ( \mu_{i},\delta_{r_{i}}) \big\} ,\ \nu= \big\{ (\nu_{j},\delta_{s_{j}}) \big\}.$$

\begin{Relation}\label{A} We have
\begin{eqnarray}
\Lan {\mu}\Big|\tau_0([D_0]\cdot \delta) \cdot\omega\Big| {\nu}\Ran_{g,\beta} C
 &=& \nonumber
\Lan {\mu}\Big|\tau_0([D_0] \cdot \delta) \cdot \omega   \prod_j 
{\tau_{\nu_{j}-1}([D_{\infty}]\cdot\delta_{s_{j}})}\Ran_{g,\beta} \\ \nonumber 
& &-\hspace{-10pt}  
\sum_{\substack{ \tilde{R}_{g,\beta} \text{ distinguished type II  }\\\tilde{R} 
\stackrel{\circ}{<} R}} 
\tilde{R} \ C_{\tilde{R},R}\\
\nonumber& & -\hspace{-10pt} \nonumber
\sum_{\substack{\|\omega^{\prime}\|\leq\|\omega\| \\
\deg(\mu^{\prime}) \geq \deg(\mu) + 1}} \hspace{-25pt}
C_{\mu^{\prime},\omega'}\Lan {\mu^{\prime}}\Big|\omega^{\prime}
\prod_j {\tau_{\nu_{j}-1}([D_{\infty}]\cdot \delta_{s_{j}})}\Ran_{g,\beta}\\ 
\nonumber\\
& & - \ \ldots \nonumber,
\end{eqnarray}
where 
$$C = \prod_j \frac{1}{(\nu_j-1)!} \ \Big( \int_{\beta} [D_\infty]\Big)^{Id(\nu)}\neq 0 $$
and the coefficients $C_{*,*}$ are fiber class integrals. The dots stand for
non-principal terms of type I and II.
\end{Relation}

\begin{proof}
Consider the Gromov-Witten invariant of the
pair $(Y,D_0)$ obtained by replacing the 
relative conditions $\nu$ along $D_{\infty}$ of $R$ by the insertions 
$$\theta=\prod_j {\tau_{\nu_{j}-1}([D_{\infty}]\cdot \delta_{s_{j}})}.$$
The resulting type I invariant, 
\begin{equation}\label{cxsd}
\Lan {\mu}\Big|\tau_0([D_0]\cdot \delta) \cdot \omega  
\prod_j {\tau_{\nu_{j}-1}([D_{\infty}]\cdot\delta_{s_{j}})}\Ran_{g,\beta},
\end{equation}
is the first term on the right side
of Relation \ref{A}. We will obtain Relation \ref{A} by computing
\eqref{cxsd} via degeneration to the normal
cone of $D_\infty$.

The special fiber of the
degeneration to the normal cone
of $D_{\infty}$ is a union of two copies of $Y$, 
$$ Y_{1}\cup_{D} Y_{2},$$ 
along a divisor $D$.
The intersection $D$ is identified with $D_{\infty}$ on $Y_{1}$ 
and $D_{0}$ on $Y_{2}$.

The degeneration
formula expresses \eqref{cxsd} in terms of type II invariants on $Y_{1}$ and 
type I invariants on $Y_{2}$ relative to $D_{0}$:
\begin{equation*}
\Lan {\mu}\Big|\tau_0([D_0]\cdot \delta) \cdot \omega  \
\theta \Ran_{g,\beta} = 
\sum
\Lan \mu 
\Big|\tau_0([D_0]\cdot \delta) \cdot\omega_{1}\Big|{\eta}
\Ran^\bullet_{{g}_{1},\mathbf{\beta}_{1}}\
 {\mathfrak{z}}(\eta)\ 
\Lan {\eta^{\vee}}\Big|\omega_{2}\ \theta
\Ran^\bullet_{{g}_{2},\mathbf{\beta}_{2}}
\end{equation*}
The sum on the right is over all 
splittings  of
$g$ and $\beta$, all distributions of the insertions of $\omega$, all 
intermediate cohomology weighted partitions ${\eta}$, and all 
configurations of connected components yielding a connected total domain.
The Gromov-Witten invariants on the right are possibly disconnected ---
indicated by the superscript $\bullet$. 
The subscript $g_i$ denotes the arithmetic  
genus of the total map to $Y_i$.

Let $\{(\eta_{k},\rho_{k})\}$
be the parts of ${\eta}$.
We may assume the insertions of $\omega$ with 
cohomology classes divisible by $[D_{0}]$ and $[D_{\infty}]$ are distributed 
to $Y_{1}$ and $Y_{2}$ respectively.
Hence, the invariants of $Y_1$ are distinguished.
For a given distribution, 
$$\int_{\beta_1}[D_{\infty}] = \int_{\beta_2}[D_{0}].$$  
Since we are omitting non-principal terms,
 we may assume either $\beta_{1}= \beta$ or $\beta_{2} = \beta$.  

\vspace{10pt}
\noindent {\bf Case 1:} $\beta_{1}=\beta$.
\vspace{10pt} 

The principal 
terms from $Y_{1}$ will be shown to be either $R$ or a distinguished invariant of type II 
lower than $R$ in our ordering.  

Let $f_i:C_i \rarr Y_i$ be the elements of the relative moduli spaces for a fixed splitting.
The condition
$\beta_{1}=\beta$ forces $\beta_{2}$ to be a multiple of the fiber class $[F]$.  Let
$\ell(\eta)$ denote the {length} of $\eta$. We find
 $$g = g_{1}+g_{2}+\ell(\eta) - 1.$$
  Since $\beta_{2}$ is a fiber class, every connected 
component of $C_2$ intersects $D_{0}$ and contains at least one relative marking. Hence, 
$$g_{2} \geq 1-\ell(\eta).$$ 
We conclude
 $g \geq g_{1}$ with equality if and only if $C_{2}$ 
consists of rational components, each totally ramified over $D_{0}$.  Since an invariant
in the degeneration formula 
with $g> g_1$ is non-principal, we consider only the extremal configurations. 

If any of the 
insertions of $\omega$ are distributed to 
$Y_{2}$, the type II invariant of $Y_{1}$ will be strictly lower than $R$
with respect to $\stackrel{\circ}{<}$. These principal contributions
appear in the second term on the right of Relation \ref{A}.

We must analyze 
the case in which $C_2$ consists of rational components totally ramified over $D_{0}$
and the only 
non-relative insertions on $Y_2$ are given by $\theta$. 

The distribution of the $\ell(\nu)$ insertions of $\theta$ among the $\ell(\eta)$ rational
components of $C_2$ decomposes ${\nu}$ into $\ell(\eta)$ 
 cohomology weighted partitions $${\nu} = \coprod_{k=1}^{\ell(\eta)} \pi^{(k)},$$
where we allow empty weighted partitions.
Here,
 $$\pi^{(k)} = \{(\nu_{n_{1}}^{(k)},\delta_{i_{1}^{(k)}}),
\dots,(\nu_{n_{s}}^{(k)},\delta_{i_{s}^{(k)}})\}.$$
Then, for each $k$,
\begin{equation}\label{zz} 
\deg({\pi^{(k)}})=\sum \deg(\delta_{i_{j}^{(k)}}) \leq  \deg(\rho_{k}).
\end{equation}
We conclude
\begin{equation*}
\deg({\nu}) \leq \deg({\eta}).
\end{equation*}  By the ordering $\stackrel{\circ}{<}$, a 
strict inequality in \eqref{zz} implies a strictly lower invariant. We 
consider only the case 
where equality holds in \eqref{zz} for each $k$.

The dimension constraint for $Y_{2}$ yields the equality
\begin{equation}\label{ddd}
\eta_k-1 = \sum_{j=1}^{\ell(\pi^{(k)})}(\pi^{(k)}_{j}-1)
\end{equation}
on each component of the domain $C_2$.

Consider the weighted partition ${\pi^{(k)}}$ containing the largest element 
$(\nu_{1},\delta_{i_{1}})$ of ${\nu}$ in the size ordering. By formula \eqref{ddd},
either $$\eta_{k} > \nu_{1}$$ or $\eta_{k} = \nu_{1}$ and all the other pairs of
$\pi^{(k)}$ are of the form $(1,\delta)$. 
In the second case, either 
$$\deg(\rho_{k}) > \deg(\delta_{i_{1}})$$ or $\rho_{k} = \delta_{i_{1}}$
and all the other pairs of $\pi^{(k)}$ are of the form $(1,\delta_{Id})$. 
Therefore,
either ${\eta}$ 
is larger than ${\nu}$ in the lexicographic ordering and corresponds to a type II invariant 
strictly lower than $R$ in the $\stackrel{\circ}{<}$ ordering, 
or the maximal pairs of $\eta$ and $\nu$ agree.  

We now repeat the above analysis for the second largest element of $\nu$ and continue 
until all the elements of $\eta$ are exhausted. We find either a strictly smaller type II
invariant in the $\stackrel{\circ}{<}$ ordering or 
 $${\eta} = {\nu}.$$ 
In the latter case, we recover $R$.

The normalization of $\theta$ sets the coefficient of $R$ in the degeneration formula to equal
$$\prod_{j} \frac{1}{(v_j-1)!} \ \Big( \int_\beta[D_\infty]\Big) ^{Id(\nu)}.$$
The coefficient is a
product of $\mathfrak{z}(\nu)$ with 
genus 0 fiber class integrals.
The fiber class integrals are evaluated using the Gromov-Witten/Hurwitz
correspondence of \cite{OP1} and the divisor equation.

\vspace{10pt}
\noindent {\bf Case 2:} $\beta_{2}=\beta$.
\vspace{10pt} 
 
The principal 
terms from $Y_{2}$ will be shown to be type I invariants of the form of the third
term on the right of Relation \ref{A}.

Let $f_i:C_i \rarr Y_i$ be the elements of the relative moduli spaces for a fixed splitting.
The condition
$\beta_{2}=\beta$ forces $\beta_{1}$ to be a multiple of the fiber class $[F]$.  
As before, after neglecting lower terms, we may assume $C_1$ consists of
$\ell(\eta)$ rational components, each totally ramified over $D_{\infty}$.

The distribution of the $\ell(\mu)$ relative markings  among the $\ell(\eta)$ rational
components of $C_1$ decomposing ${\mu}$ into $\ell(\eta)$ 
 cohomology weighted partitions $${\mu} = \coprod_{k=1}^{\ell(\eta)} \pi^{(k)},$$
where empty weighted partitions are {\em not} allowed.
 
If the $k$th  component of $C_1$ does not
contain the distinguished marked point, then
$$\deg({\pi^{(k)}}) + \deg(\rho_{k}) \leq \text{dim}_{\mathbb R} (X)$$ 
since all these classes are pulled-back from the same projection map to $X$.
If the $k$th component of $C_1$ does
contain the distinguished marked point, then
$$\deg({\pi^{(k)}}) + \deg(\rho_{k}) \leq \text{dim}_{\mathbb R}(X) -1$$ 
since there is an additional class of nonzero degree from the distinguished marking.
The result 
follows 
since the cohomology weights of ${\eta^{\vee}}$ are Poincar\'e dual to the classes $\rho_{k}$.
\end{proof}

\subsection{Rubber calculus}\label{rc}
\subsubsection{Rubber targets}
We will  study the Gromov-Witten theory of the pair $(Y,D_0)$
via virtual localization
\cite{GP,GV} with respect to the natural fiberwise  $\com^{\ast}$-action  on $Y$ discussed in
Section \ref{fca}.  

The $\com^*$-action on $Y$ induces a canonical $\com^*$-action on the moduli space of
stable relative maps to the pair $(Y,D_0)$.
The
 $\com^*$-fixed loci of the latter action 
involve stable relative maps to non-rigid targets.  
Let $$\overline{M}^\bullet=
\overline{M}^\bullet_{g,n}(Y/D_{0}\cup D_{\infty},\beta,\stackrel{\rarr}{\mu},\stackrel{\rarr}{\nu})$$
denote the moduli space of stable maps 
to $Y$ relative to both divisors, and
let 
$$\overline{M}^{\bullet\sim}=\overline{M}_{g,n}^{\bullet \sim}(Y/D_{0}\cup D_{\infty},\beta,
\stackrel{\rarr}{\mu},\stackrel{\rarr}{\nu})$$ 
denote the corresponding space of stable maps to a non-rigid target ---
termed a {\em rubber} target in \cite{FP2,OP3}.
Let 
$$\epsilon:\overline{M}^\bullet \rightarrow \overline{M}^{\bullet \sim}$$
be the canonical forgetful map.

The superscripted $\bullet$ indicates
the moduli spaces $\overline{M}^\bullet$ and $\overline{M}^{\bullet \sim}$  may
parameterize maps with disconnected domains  with specified genus and class distributions.
The latter data is not made explicit in our notation. The subscripted $g$
is the arithmetic genus of the total domain.
Similarly, the brackets  $\lan,\ran^\bullet$ and $\lan,\ran^{\bullet \sim}$ 
will denote invariants with possibly disconnected domains.
There is no product rule relating connected and disconnected rubber invariants.

\subsubsection{Cotangent classes}
The moduli space $\overline{M}^\bullet$ and $\overline{M}^{\bullet\sim}$
carry tautological cotangent line bundles ${\mathbb{L}}_0$ and ${\mathbb{L}}_\infty$
determined by the relative divisors. The associated cotangent line classes
$$\Psi_0= c_1({\mathbb{L}}_0), \ \ \Psi_\infty = c_1({\mathbb{L}}_\infty),$$
play an important role in relative Gromov-Witten theory.

 Let ${\text{pr}}_1$ and ${\text{pr}}_2$
denote the projections onto the first and second factors of
the product
$$D_0 \times \overline{M}^\bullet.$$
Let 
$\tau: T \rarr \overline{M}^\bullet$
denote the universal family of {\em targets} over the moduli space, and let
$$\iota: D_0 \times \overline{M}^\bullet \rarr T$$
denote the inclusion of the relative divisor.
The cotangent line determined by $D_0$ is defined by
\begin{equation}\label{sss12}
{\mathbb L}_0 = {\text{pr}}_{2*}\Big({\text{Conorm}}(\iota) \otimes{\text{pr}_1}^*(\text{Norm}(Y/D_0))
\Big),
\end{equation} 
 where $\text{Conorm}(\iota)$ is the 
conormal bundle of the embedding $\iota$ and
$$\text{Norm}(Y/D_0) = L^*$$
is the normal bundle of $D_0$ in $Y$.
The push-forward \eqref{sss12} is easily seen to define a {\em line bundle} ${\mathbb L}_0$.

The line bundle ${\mathbb L}_\infty$ on $\overline{M}^\bullet$
is similarly defined. The constructions in the rubber case are
identical.

\subsubsection{Rigidification}
The following rigidification Lemma plays a fundamental role in our localization analysis.

\begin{Lemma}  \label{ss12}
Let $p$ be 
a non-relative marking with evaluation map $$\mathrm{ev}_p:\overline{M}^\bullet \rightarrow Y.$$ 
Then,
\begin{eqnarray*}
[\overline{M}^{\bullet \sim}]^{\mathrm{vir}} & =&  \epsilon_{\ast}\Big(\ev^{\ast}_p([D_{0}])\cap 
[\overline{M}^\bullet]^{\mathrm{vir}}\Big) \\
& =& \epsilon_{\ast}\Big(\ev^{\ast}_p([D_{\infty}])\cap [\overline{M}^\bullet]^{\mathrm{vir}}\Big).
\end{eqnarray*}
\end{Lemma}
\begin{proof}
The forgetful map $\epsilon$ is equivariant with respect to the canonical
$\com^*$-action on $\overline{M}^\bullet$ induced from the fiberwise $\com^*$-action on $Y$ and 
the trivial
 $\com^*$-action on $\overline{M}^{\bullet\sim}$.
We prove the first equality by $\com^*$-localization.

 A stable 
relative map corresponding to an element of
 a typical $\com^*$-fixed locus of $\overline{M}^\bullet$ is a union of three basic
submaps:
\begin{enumerate}
\item[(i)]a 
nonrigid stable
map to the degeneration of $Y$ over $D_0$,
\item[(ii)] a nonrigid stable map to the degeneration of $Y$ over $D_\infty$,
\item[(iii)]
a collection of $\com^*$-invariant, fiber class, rational Galois covers joining (i-ii).
\end{enumerate}
The forgetful map simply contracts the intermediate 
rational curves (iii). 

Assuming we have a proper degeneration on each side of $Y$, the virtual dimension of the 
$\com^*$-fixed locus is
2 less than the virtual dimension of $\overline{M}^\bullet$.
Since the dimension of 
\begin{equation}\label{cgt}
\epsilon_{\ast}\Big(\ev^{\ast}([D_{0}])\cap [\overline{M}^\bullet]^{\mathrm{vir}}\Big)
\end{equation}
is only 1 less than the virtual dimension of $\overline{M}^\bullet$, the contributions of the above 
loci cancel in the computation of the  push-forward \eqref{cgt}.

The only fixed loci which may contribute are
those with target degeneration on only one side of $Y$.
Since the 
marking $p$ is constrained by an insertion of $[D_{0}]$, we need only consider degenerations
along $D_{0}$.  

There is a 
unique $\com^*$-fixed locus which provides non-cancelling contributions to the push-forward
\eqref{cgt}. Moreover, the $\com^*$-fixed
locus is isomorphic to $\overline{M}^{\bullet\sim}$. The contribution,
$$\frac{-\ev_p^*(c_1(L))+t}{-\Psi_\infty+t}\cap [\overline{M}^{\bullet\sim}]^{vir},$$
is obtained from the virtual localization formula.
Here, $\Psi_\infty$ is the cotangent line class on $\overline{M}^{\bullet \sim}$ at the relative
divisor $D_\infty$. By dimension considerations, the only non-cancelling part
is $[\overline{M}^{\bullet \sim}]^{vir}$ --- proving 
the first equality. The proof of the second equality is identical.
\end{proof}
\label{rbb}
\subsubsection{Dilation and divisor}
As before,
let $\Psi_\infty$ denote the cotangent line class on $\overline{M}^{\bullet \sim}$ at the relative
divisor $D_\infty$. 
The dilaton equation for rubber integrals is
\begin{multline*}
\Lan \mu \Big| \tau_1(1) \prod_{i=1}^n\tau_{k_i}(\gamma_{l_i}) 
\Psi_\infty^k \Big| \nu \Ran_{g,\beta}^{\bullet \sim} 
= \\
\Big(2g-2+n+\ell(\mu)+\ell(\nu)\Big) 
\Lan \mu \Big| \tau_1(1) 
\prod_{i=1}^n\tau_{k_i}(\gamma_{l_i}) \Psi_\infty^k \Big| \nu \Ran_{g,\beta}^{\bullet \sim}.
\end{multline*}
The divisor equation for $H\in H^2(X,\Q)$, however, takes a modified form:
\begin{multline*}
\Lan \mu \Big| 
\tau_0(H) \prod_{i=1}^n\tau_{k_i}(\gamma_{l_i}) \Psi_\infty^k \Big| \nu \Ran_{g,\beta}^{\bullet \sim} 
=\\
\Big(\int_{\pi_*(\beta)} H\Big) \cdot
\Lan \mu \Big| 
\prod_{i=1}^n\tau_{k_i}(\gamma_{l_i}) \Psi_\infty^k \Big| \nu \Ran_{g,\beta}^{\bullet \sim} \\ 
\ \ \ \ \ \  + \sum_{j=1}^n \Lan \mu \Big| \ldots \tau_{k_j-1}(\gamma_{l_j} \cdot H) \ldots 
\Psi_\infty^k \Big| \nu \Ran_{g,\beta}^{\bullet \sim} \\ 
+ \sum_{j=1}^{\ell(\nu)} 
\Lan \mu \Big| \prod_{i=1}^n\tau_{k_i}(\gamma_{l_i}) \Psi_\infty^{k-1} \Big| \{ 
\ldots (\nu_j,\delta_{s_j}\cdot H) \ldots \}\Ran_{g,\beta}^{\bullet \sim}\cdot  \nu_j
\end{multline*}
The dilation and divisor equations are proven by the standard cotangent
line comparison method.

\subsubsection{Calculus I: fiber class}
The rubber calculus relates Gromov-Witten rubber invariants with
$\Psi_\infty$ insertions to Gromov-Witten
invariants of the pair $(Y, D_0 \cup D_\infty)$.

Consider first a rubber integral with descendent insertions $\omega$
for which $\beta$ is a multiple of the
fiber class:
\begin{equation}
\label{xxbb}
\Lan \mu \Big|  \omega\ \Psi_\infty^k \Big| \nu \Ran_{g,\beta}^{\bullet \sim}.
\end{equation}
A contracted genus 0 component of the domain  must carry at least 3 non-relative markings
by stability. 
Similarly, a contracted genus 1 domain component must carry at least 1 non-relative marking.
A non-contracted domain component  must carry at least 2 relative markings ---
the intersection points with  $D_0$ and $D_\infty$. 
Finally, by target stability, not all domain components can be genus 0 and fully ramified over 
$D_0$ and $D_\infty$.
We conclude
$$2g-2+n+\ell(\mu)+\ell(\nu)>0.$$
Therefore, the fiber class rubber integral \eqref{xxbb} is determined
by
\begin{equation}\label{gfww1}
\Lan 
\mu \Big| \tau_1(1) \cdot \omega\ \Psi_\infty^k \Big| \nu \Ran_{g,\beta}^{\bullet \sim}.
\end{equation}
and the dilaton equation.

Let $p$ denote 
the marked point carrying the insertion $\tau_1(1)$ in the rubber integral \eqref{gfww1}.
There is canonical map
to the Artin stack of genus 0, 3-pointed curves,
$$\alpha:\overline{M}^{\bullet \sim} \rarr \mathcal{M}_{0,3}.$$
Given $[f]\in \overline{M}^\sim$, $\alpha(f)$ is the genus 0 curve 
$$C_f=\pi^{-1}\big(\pi(f(p))\big)$$
with 3 marking determined by
$$D_0 \cap C_f,\ f(p),\ D_\infty \cap C_f.$$
The class $$\Psi_\infty- \text{ev}_p^*(c_1(L))$$ is the pull-back of the cotangent line of the third
marking on the Artin stack.

Topological recursion relation with respect to the cotangent class of the third
marking of  $\mathcal{M}_{0,3}$  can be 
pulled-back via $\alpha$:
\begin{multline*}
\Lan 
\mu \Big| \tau_1(1) \cdot \omega\ \Psi_\infty^k \Big| \nu \Ran_{g,\beta}^{\bullet \sim} =
\Lan 
\mu \Big| \tau_1(c_1(L)) \cdot \omega\ \Psi_\infty^{k-1} \Big| \nu \Ran_{g,\beta}^{\bullet \sim}
\\
+ \sum\Lan {\mu}\Big|\tau_1(1)\cdot \omega_1
\Big| {\eta}\Ran^{\bullet \sim}_{g_{1},\beta_{1}} 
\mathfrak{z}(\eta)\Lan \eta^{\vee}\Big| \omega_2 \
\Psi_{\infty}^{k-1}\Big|\nu\Ran^{\bullet \sim}_{g_{2},\beta_{2}},
\end{multline*}
where the sum is over all splitting of $g$ and $\beta$, all distributions of the insertions, and
all intermediate cohomology weighted partitions.

The first term of the sum on the right can be expressed as a type II
invariant by Lemma \ref{ss12}:
$$
\Lan {\mu}\Big|\tau_1(1)\cdot \omega_1
\Big| {\eta}\Ran^{\bullet \sim}_{g_{1},\beta_{1}} =
\Lan {\mu}\Big|\tau_1([D_0])\cdot \omega_1
\Big| {\eta}\Ran_{g_{1},\beta_{1}}^\bullet.$$
For the application of Lemma \ref{ss12} here, we require the compatibility of $\epsilon$ with the
cotangent lines at the marked points. 

We have reduced the original fiber class rubber invariant \eqref{xxbb} to 
invariants of the same type with {\em fewer}  $\Psi_\infty$ insertions.
Repeating the cycle yields rubber invariants without $\Psi_\infty$ insertions.
The latter are related to type II invariants by Lemma \ref{ss12} after adding a dilaton
insertion.

\subsubsection{Calculus II: $\pi_*(\beta)\neq 0$}
If $\beta$ is not fiber class, then $\pi_*(\beta)\neq 0$.
Consider a non-fiber class rubber integral:
\begin{equation}
\label{xxxbb}
\Lan \mu \Big|  \omega\ \Psi_\infty^k \Big| \nu \Ran_{g,\beta}^{\bullet \sim}.
\end{equation}
If $H\in H^2(X,\Q)$ is an ample class, then
$$\int_{\pi_*(\beta)}H >0.$$
Therefore, by the divisor equation, the rubber integral
\begin{equation}\label{gfww2}
\Lan 
\mu \Big| 
\tau_0(H) \cdot \omega \ \Psi_\infty^k \Big| \nu \Ran_{g,\beta}^{\bullet \sim}.
\end{equation}
determines \eqref{xxxbb} modulo rubber integrals with strictly fewer cotangent lines.

As in the fiber case, we may apply the topological recursion relations to 
\eqref{gfww2}. By repeating the cycle and applying Lemma \ref{ss12} after all $\Psi_\infty$
insertions are removed, we can express the original rubber invariant \eqref{xxxbb}
in terms of type II invariants.

A refined consequence of the rubber calculus will be needed in the 
proof of Relation \ref{B} in the following 
Section.
\begin{Lemma}\label{ab}
A rubber invariant with {\em connected} domain
and class satisfying $\pi_*(\beta)\neq 0$
is expressed by the calculus in terms of 
type II invariants as:
\begin{equation}
\label{xfw}
 \Lan \mu \Big|  \omega\ \Psi_\infty^k \Big| \nu \Ran_{g,\beta}^\sim = 
\sum_{\substack{
\|\omega^{\prime}\| \leq \|\omega\|\\
\deg(\mu^{\prime})\geq \deg(\mu)\\
\deg(\nu^{\prime}) \geq \deg(\nu)\\ m\geq 0}} 
C_{\mu,'\omega',\nu'}\Lan \mu' \Big|  \tau_0([D_0]\cdot H \cdot c_1(L)^{m}) \cdot \omega'  \Big|
 \nu' \Ran_{g,\beta} + \ \ldots\ .
\end{equation}
The brackets $\lan,\ran$ on the right denote {\em connected} invariants.
The coefficients $C_{*,*,*}$ are determined by fiber class integrals, 
and the dots stand for non-principal terms of type II.
\end{Lemma}

\subsection{Relation \ref{B}}

The next relation expresses the 
type I invariants occurring in Relation \ref{A} in terms of  type II invariants and
the Gromov-Witten theory of $X$.

Consider the type I invariant 
$$\Lan{\mu}\Big| \tau_0([D_0]\cdot \delta)\cdot 
\omega \prod_j {\tau_{\nu_{j}-1}([D_{\infty}]\cdot 
\delta_{s_{j}})}\Ran_{g,\beta}$$
with $\pi_*(\beta)\neq 0$, $\deg(\delta)>0$, and 
relative conditions
$$\mu=\big\{ ( \mu_{i}, \delta_{r_{i}}) \big\}, \ \nu= \big\{ (\nu_{j},\delta_{s_{j}}) \big\}.$$
\begin{Relation}\label{B} We have,
\begin{multline*} 
\Lan{\mu}\Big| \tau_0([D_0]\cdot \delta)\cdot 
\omega \prod_j {\tau_{\nu_{j}-1}([D_{\infty}]\cdot 
\delta_{s_{j}})}\Ran_{g,\beta} =
\\  
\ \ \  \sum_{\substack{\|\omega^{\prime} \|\leq \|\omega\|\\
\deg(\mu^{\prime})\geq \deg(\mu)+1 \\
\deg(\nu^{\prime}) \geq \deg(\nu)\\m\geq 0}} 
C_{\mu^{\prime},\nu^{\prime},
\omega^{\prime}} \Lan {\mu^{\prime}}\Big|\tau_0([D_{0}]\cdot H \cdot c_1(L)^m)\cdot \omega^{\prime}\Big|
{\nu^{\prime}}\Ran_{g,\beta}
\\  
+\sum_{\substack{\|\omega^{\prime} \|\leq \|\omega\|\\
\deg(\mu^{\prime})\geq \deg(\mu) \\
\deg(\nu^{\prime}) \geq \deg(\nu)+1\\ m\geq 0}} 
C_{\mu^{\prime},\nu^{\prime},
\omega^{\prime}} \Lan {\mu^{\prime}}\Big|\tau_0([D_{0}]\cdot H\cdot c_1(L)^{m})\cdot \omega^{\prime}\Big|
{\nu^{\prime}}\Ran_{g,\beta}
\\ \ \ \ \ \ \ \ \ \ \ \ \ \ \ 
-\sum_{\substack{\|\omega^{\prime} \|\leq \|\omega\|\\
\deg(\mu^{\prime})\geq \deg(\mu) \\
\deg(\nu^{\prime}) \geq \deg(\nu)\\ m\geq 0}} 
C_{\mu^{\prime},\nu^{\prime},
\omega^{\prime}} 
\Lan {\mu^{\prime}}\Big|\tau_0([D_{0}]\cdot c_1(L)^{m+1}\cdot \delta)\cdot \omega^{\prime}\Big|
{\nu^{\prime}}\Ran_{g,\beta}
+\ \ldots.\\
\end{multline*}
where $H\in H^2(X,\Q)$ is an ample class and the 
coefficients $C_{*,*,*}$ are fiber class integrals. 
The dots stand for non-principal terms of type II and integrals in the Gromov-Witten theory of $X$
\end{Relation}

\begin{proof}
The first step is to use the basic divisor relation in $H^2(Y,\Q)$,
$$[D_0]=[D_\infty]-c_1(L),$$
to rewrite the invariant on the left side of Relation \ref{B}:
\begin{multline}\label{dcq}
\Lan{\mu}\Big| \tau_0([D_0]\cdot \delta)\cdot 
\omega \prod_j {\tau_{\nu_{j}-1}([D_{\infty}]\cdot 
\delta_{s_{j}})}\Ran_{g,\beta} = \\ 
\ \ \ \ \ \ \ \ \ \ \  \ \ \ \ \ \  \ \ \ \ \ \ \ \ \ \ \ \ \ \ \ \ \  \ \  
\Lan{\mu}\Big| \tau_0([D_\infty]\cdot \delta)\cdot 
\omega \prod_j {\tau_{\nu_{j}-1}([D_{\infty}]\cdot 
\delta_{s_{j}})}\Ran_{g,\beta} \\ \nonumber
- \Lan{\mu}\Big| \tau_0(c_1(L)\cdot \delta)\cdot 
\omega \prod_j {\tau_{\nu_{j}-1}([D_{\infty}]\cdot 
\delta_{s_{j}})}\Ran_{g,\beta}.
\end{multline}
We will calculate the two latter type I invariants  by localization
on $(Y,D_0)$.

The fiberwise $\com^*$-action on $Y$ induces a canonical action on the
moduli space of stable relative maps to the pair $(Y,D_0)$.
The $\com^*$-fixed loci 
consist of  moduli space of stable maps to rubber over $D_0$ connected by
fiberwise rational Galois covers to  moduli spaces of
stable maps to $D_\infty$. The connection data of the Galois covers is
described by a sum over cohomology weighted partitions $\eta$ 
specifying the rubber relative conditions on the connecting divisor.

We first study the localization calculation of the type I invariant
\begin{equation}\label{fz34}
\Lan{\mu}\Big| \tau_0([D_\infty]\cdot \delta)\cdot 
\omega \prod_j {\tau_{\nu_{j}-1}([D_{\infty}]\cdot 
\delta_{s_{j}})}\Ran_{g,\beta}.
\end{equation}
The insertions of \eqref{fz34} all have canonical equivariant lifts.
The
marked points not included in $\omega$ are all distributed to the moduli space
of stable maps to $D_\infty$.  
By the localization formula,
the contribution of the latter moduli space is
simply a Hodge integral in the Gromov-Witten theory of $X$. The Hodge class
may be removed by \cite{FP}.

Since we neglect lower terms, we only need to consider $\com^*$-fixed loci 
for which the maps to rubber are of genus $g$ and class $\beta$.
Let $C_0$ be the subcurve of the domain mapping to rubber, and let
$C_\infty$ be the subcurve mapping to $D_\infty$.
Since all components of $C_\infty$ are contracted, 
the argument used for rational fibers in the proof of Relation \ref{A} shows
 $$\deg({\eta}) \geq \deg({\nu})+\deg(\delta).$$
The principal terms of the localization formula for \eqref{fz34} are
therefore of the form
$$\Lan \mu \Big| \omega'\ \Psi_\infty^k \Big| \eta \Ran_{g,\beta}^\sim$$
where $\| \omega' \| \leq \|\omega \|$ and 
$\deg(\eta) \geq \deg(\nu)+1$.

By Lemma \ref{ab},  the invariant \eqref{fz34} contributes only principal terms
of the type of the second summand on the right side of Relation \ref{B}.

Next, we study the the localization calculation of the type I invariant
\begin{equation}\label{xcf12}
\Lan{\mu}\Big| \tau_0(c_1(L)\cdot \delta)\cdot 
\omega \prod_j {\tau_{\nu_{j}-1}([D_{\infty}]\cdot 
\delta_{s_{j}})}\Ran_{g,\beta}.
\end{equation}
If the insertion $\tau_0(c_1(L)\cdot\delta)$ is distributed to $D_\infty$, then,
as above, we
only obtain principal terms of the type of the second summand on the
right  of Relation \ref{B}.

If the insertion $\tau_0(c_1(L)\cdot\delta)$ is distributed to $D_0$, then,
by the rubber calculus, we 
only obtain principal terms of the type of the first and third summands on the
right  of Relation \ref{B}.
\end{proof}

We will also require a version of Relation \ref{B} {\em without} the 
distinguished insertion $\tau_0([D_0] \cdot \delta)$. The proof is
identical.

\vspace{+10pt}
\noindent {\bf Relation ${\mathbf 2'}$.}
{\em We
 have}
\begin{multline*} 
\Lan{\mu}\Big| 
\omega \prod_j {\tau_{\nu_{j}-1}([D_{\infty}]\cdot 
\delta_{s_{j}})}\Ran_{g,\beta} =
\\  
\ \ \  \sum_{\substack{\|\omega^{\prime} \|\leq \|\omega\|\\
\deg(\mu^{\prime})\geq \deg(\mu) \\
\deg(\nu^{\prime}) \geq \deg(\nu)\\ m \geq 0}} 
C_{\mu^{\prime},\nu^{\prime},
\omega^{\prime}} \Lan {\mu^{\prime}}\Big| \tau_0([D_0]\cdot H\cdot c_1(L)^m)\cdot  \omega^{\prime}\Big|
{\nu^{\prime}}\Ran_{g,\beta} +\ldots ,
\end{multline*}
{\em where
the dots stand for non-principal terms of
type II and integrals in the Gromov-Witten theory of $X$.}

\subsection{Proof of Theorem \ref{mn}} 
Our primary induction is on the pair $(g,\beta)$ where
$$(g',\beta')<(g,\beta)$$
if $\beta'<\beta$ or if $\beta'=\beta$ and $g'<g$.
If $\pi_*(\beta)=0$, the invariants are fiber class and are determined for all $g$
by Section \ref{fclass}.

By a secondary induction on the $\stackrel{\circ}{<}$ ordering, 
Relations 1, 2, and $2'$ determine all distinguished invariants of type II.
Consider an invariant
\begin{equation}\label{xxp2}
\Lan {\mu}\Big|\tau_0([D_0]\cdot \delta) \cdot\omega\Big| {\nu}\Ran_{g,\beta}
\end{equation}
for which $\pi_*(\beta) \neq 0$.
We apply Relation 1 to \eqref{xxp2}.
To the first term on the right, we apply Relation 2. To the third term on
the right, we apply relation $2'$. The outcome modulo the primary induction
 is a determination of 
\eqref{xxp2} in terms of 
distinguished type II invariants lower in the $\stackrel{\circ}{<}$ ordering.
The proof of Theorem \ref{mn} for distinguished type II invariants is complete.

By localization and the rubber calculus, every type I and II invariant can be expressed
in terms of distinguished type II invariants and Hodge integrals in the
Gromov-Witten theory of X.
\qed

\subsection{Connected/disconnected invariants} \label{disco}
For simplicity, we 
have assumed the disconnected invariants of the pairs $(Y,D_0)$,
$(Y,D_\infty)$, and $(Y,D_0\cup D_\infty)$
factor as a product of connected invariants.
Unfortunately, the product rule for disconnected invariants has not been 
included in the foundational treatments of the subject. 

The entire proof of Theorem \ref{mn} is valid {\em without} assuming
the product rule. The induction argument reduces the disconnected 
invariants of the three relative pairs to the Gromov-Witten theory of
$X$. 

If the product rule is not assumed, 
the only difference in the proof of Theorem \ref{mn}
occurs in the treatment of the fiber class integrals in
Section \ref{fclass}. Disconnected fiber class invariants
must be considered.
Let 
$$P=\bigcup_{i=1}^n \proj^1$$
be a disconnected set of projective lines.
Denote the disconnected divisors 
$$\{0,\ldots,0\}, \ \{\infty, \ldots, \infty\} \subset P$$
by $D_0$ and $D_\infty$.
The disconnected fiber class invariants require the computation of the 
$(\com^*)^n$-equivariant Gromov-Witten theory of the
pairs $$(P,D_0),\ (P,D_\infty),\ (P,D_0\cup D_\infty).$$
An elementary $(\com^*)^n$-localization argument reduces the 
study of $P$ to the product of $n$ copies of $\proj^1$ ---
as predicted by the product rule. We leave the details to
the reader.
 
Indeed, our study of relative Gromov-Witten theory can be used
to {\em prove} the product rule. The argument will be presented
elsewhere.

\section{Relative in terms of absolute} \label{bbb}
\subsection{Notation}
Let $V$ be a nonsingular,  
complex, projective variety containing a nonsingular divisor $W$,
$$\iota:W \rarr V.$$
Let $\iota^*$ denote the restriction map on cohomology,
$$\iota^*: H^*(V,\Q) \rarr H^*(W,\Q).$$
The cohomological push-forward 
$$\iota_*: H^*(W,\Q) \rarr H^*(V,\Q)$$ is
determined by the restriction map $\iota^*$ and
Poincar\'e duality.

Let $N$ be the normal bundle of $W$ in $V$. Since
$$c_1(N)= \iota^*(c_1(T_V))-c_1(T_W)\in H^*(W,\Q),$$ 
the Chern class $c_1(N)$ is also is determined by the restriction 
map $\iota^*$.

We denote the Gromov-Witten  invariants of the pair $(V,W)$ by the right bracket
\begin{equation} \label{x123}
\left. \Lan  \prod_i \tau_{k_i}(\gamma_{l_i}) \ \right| {\mathbf \nu} 
\Ran^{V/W}_{g,\beta}
\end{equation}
where 
$$\{\gamma_1, \ldots, \gamma_{m_V}\}, \ \{\delta_1, \ldots, \delta_{m_W}\}$$
are bases of $H^*(V,\Q)$ and $H^*(W,\Q)$,
$$\mathbf{\nu}=
\left\{ (\nu_1, \delta_{s_1}), \ldots, (\nu_{\ell(\nu)}, \delta_{s_{\ell(\nu)}}) 
\right\}, $$
and $\sum_j \nu_j=\int_\beta [W]$.

The above bracket denotes Gromov-Witten invariants with
 connected domain curves. As before, we treat disconnected
invariants as products of connected invariants. Our proof of Theorem \ref{wd} is valid
without assuming 
the product rule --- see the discussion of Section \ref{disco}.

\subsection{Degeneration}

Let $\mathcal{F}$ be the degeneration to the normal cone of $W$.
The degeneration formula \cite{EGH,IP,LR,L} applied to
$\mathcal{F}$ expresses the absolute Gromov-Witten theory of
$V$ in terms of the relative theories of 
$(V,W)$ and $(\proj(N\oplus \oh_W), W)$.

By Theorem \ref{mn}, the
relative theory of $(\proj(N\oplus \mathcal{O}_{W}),W)$ is determined by
the absolute theory of $W$ and the Chern class
$c_1(N)\in H^2(W,\Q)$.
In order to prove Theorem \ref{wd}, we view the degeneration formula
as providing equations for the invariants \eqref{x123} in terms of the 
Gromov-Witten theories of $V$ and 
$(\proj(N\oplus \mathcal{O}_{W}),W)$.

\subsection{Ordering}
We partially order the invariants of $(V,W)$ by a relation very similar to
the ordering of distinguished type II invariants of Section \ref{xxw}:
$$\Lan 
 \omega' \Big |\nu'\Ran_{g',\beta'} \stackrel{\circ}{<} 
\Lan  \omega
\Big|\nu \Ran_{g,\beta}$$
if one of the conditions below holds
\begin{enumerate}
\item[(1)] $\beta' < \beta$, 
\item[(2)] equality in (1) and $g' < g$,
\item[(3)] equality in (1-2) and $\| \omega'\| < \| \omega \|$,
\item[(4)] equality in (1-3) and $\deg({\nu}') > \deg({\nu})$,
\item[(5)] equality in (1-4) and ${\nu}' \stackrel{l}{>} {\nu}$,
\end{enumerate}
Here, $\omega'$ and $\omega$ represent products of
descendent insertions.
For any given invariant of $(V,W)$ , there are only finitely many 
invariants of $(V,W)$ lower in the partial ordering.  

The degeneration equations for the invariants $(V,W)$ will be proven to
be lower triangular with respect to the $\stackrel{\circ}{<}$ ordering ---
and therefore nonsingular.

\subsection{Proof of Theorem \ref{wd}}

To each relative invariant \eqref{x123} of the pair $(V,W)$, we associate the following
 absolute
 invariant of $V$
\begin{equation} \label{x1234}
 \Lan  \prod_i \tau_{k_i}(\gamma_{l_i}) \cdot \prod_j\tau_{\nu_j-1}\big(
\iota_{*}(\delta_{s_j})\big)
\Ran^{V}_{g,\beta}.
\end{equation}

In order to evaluate the absolute invariant \eqref{x1234} by the degeneration formula,
we must lift the cohomology classes $\gamma_{l_i}$ and $\iota_*(\delta_{s_j})$ to
the total space of the family ${\mathcal F}$:
\begin{enumerate}
\item[(i)]
the classes $\gamma_{l_i}$ are lifted by pull-back via the first factor
of the blow-down map
$${\mathcal F} \rarr V\times \com,$$
\item[(ii)] the class $\iota_*(\delta_{s_j})$ are lifted by
$$\iota_{W\times\com,*}(\delta_{s_j}\otimes Id)\in H^*({\mathcal F},\Q)$$
where
$$\iota_{W\times \com}: W\times \com \rarr {\mathcal F}$$
is the inclusion in the blow-up via strict transform.
\end{enumerate}

\begin{Lemma}\label{nwq}
The principal terms of the degeneration equation are
\begin{multline*}
\Lan  \prod_i \tau_{k_i}(\gamma_{l_i}) \cdot \prod_j\tau_{\nu_j-1}\big(
\iota_{*}(\delta_{s_j})\big)
\Ran^{V}_{g,\beta} = 
\\ 
  \Lan \left. \prod_i \tau_{k_i}(\gamma_{l_i}) \  \right| {\mathbf \nu} 
\Ran^{V/W}_{g,\beta} C
+ \left. \sum
  \Lan \omega' \right| \nu'\Ran_{g,\beta}^{V/W} C_{\omega',\nu'} +\ldots,
\end{multline*}
where 
$$C= \prod_j\frac{1}{(\nu_j-1)!}\ \Big(\int_\beta[W]\Big)^{Id(\nu)}\neq 0,$$
the coefficients $C_{*,*}$ are
fiber class integrals, and the 
 sum  is over lower invariants in the
$\stackrel{\circ}{<}$ ordering.
\end{Lemma}

Following the notation of Section \ref{ax1}, the {\em principal terms} of the
equation are the invariants of $(V,W)$ of genus $g$ and class $\beta$.

\begin{proof}
The strategy is identical to the proof of Relation \ref{A}.
We need only consider degeneration splittings for which the $(V,W)$ side
carries genus $g$, class $\beta$, and {\em all} the insertions
 $$\prod_i \tau_{k_i}(\gamma_{l_i}).$$
By the lifting choice, all the insertions
$$\prod_j\tau_{\nu_j-1}\big(
\iota_{*}(\delta_{s_j})\big)$$
are distributed to the $(\proj(N\oplus \mathcal{O}_{W}),W)$ side.
The analysis of the rational fiber class integrals then
exactly follows Case 1 of the proof of Relation \ref{A}.
\end{proof}

By Lemma \ref{nwq}, the degeneration equations form a lower triangular
system determining the invariants of $(V,W)$ in terms of
the invariants of $V$ and $(\proj(N\oplus \mathcal{O}_{W}),W)$.
The proof of Theorem \ref{wd} is complete. \qed

\subsection{Proof of Corollary \ref{rwd}}
Let $(V,W)$ be a hypersurface pair in $\proj^r$.
We first prove 
$$\iota_*(\delta)\in H^*(V,\Q)$$ is
simple for any $\delta \in H^*(W,\Q)$.
If $\delta$ is simple, the result is clear. If
$\delta$ is  not simple, then
 $$\deg(\delta)=r-2$$ by Lefchetz. We conclude 
$$\deg(\iota_*(\delta))=r,$$ and thus $\iota_*(\delta)$
is simple by Lefchetz.

Consider the degeneration equation of Lemma \ref{nwq}  associated to a simple
invariant of $(V,W)$. The left side involves a simple invariant
of $V$ since $\iota_*(\delta)$ is always simple. The right side
involves simple invariants of $(V,W)$ and 
general invariants
of $(\proj(N\oplus \mathcal{O}_{W}),W)$.

The degeneration equations form a lower triangular
system determining the simple invariants of $(V,W)$ in terms of
the simple invariants of $V$ and general invariants of
$(\proj(N\oplus \mathcal{O}_{W}),W)$.
\qed

\section{Mayer-Vietoris and the quintic scheme}
\label{edd}
\subsection{Proof of Lemma \ref{ss}}
Our calculation scheme for hypersurfaces depends upon Theorem \ref{mv} and
Lemma \ref{ss}.
Theorem \ref{mv} is an immediate consequence of the degeneration formula and Theorem 
\ref{wd}.

Let $V$ be a nonsingular, projective variety. Let
$Z\subset V$ be the nonsingular complete intersection of two nonsingular divisors
$$W_1,W_2\subset V,$$
and let $\widetilde{V}$ be the blow-up of $V$ along $Z$.

\setcounter{Lemma}{0}
\begin{Lemma}\nonumber
The Gromov-Witten theory of $\widetilde{V}$ is
uniquely and effectively determined  by
the Gromov-Witten theories of $V$, $W_1$, and $Z$
and the restriction maps
$$H^*(V,\Q) \rarr H^*(W_1,\Q) \rarr H^*(Z,\Q).$$ 
\end{Lemma}

\begin{proof}
Let $N_i$ be the normal bundle to $W_i$ in $V$.
Let ${\mathcal F}$ be the degeneration 
to the normal cone of $W_1$ in $V$.
By strict transform, there is an inclusion
$$\iota_{W_1\times \com}: W_1 \times \com \rarr {\mathcal F}.$$
Let $\widetilde{\mathcal F}$ be the blow-up of ${\mathcal F}$
along $Z\times \com$ embedded via 
$$Z\times \com \subset W_1\times \com\stackrel{\iota_{W_1\times \com}}{\longrightarrow} 
{\mathcal F}.$$

The family ${\mathcal F}$ is a degeneration of ${V}$ to
$(V,W_1)$ and $({\proj}(N_1\oplus \oh_{W_1}), D_0)$.
The family $\widetilde{\mathcal F}$ is a degeneration of $\widetilde{V}$ to
$(V,W_1)$ and 
$(\widetilde{\proj}(N_1\oplus \oh_{W_1}), D_0)$ where 
$\widetilde{\proj}(N_1\oplus \oh_{W_1})$
 is the blow-up
of ${\proj}(N_1\oplus \oh_{W_1})$ along $Z\subset D_\infty$.

By the degeneration formula, the Gromov-Witten theory
of $\widetilde{V}$ is determined by the Gromov-Witten theories
of $(V,W_1)$ and $(\widetilde{\proj}(N_1\oplus \oh_{W_1}), D_0)$
since the non-vanishing cohomology for the family $\widetilde{\mathcal F}$
is all of $H^*(\widetilde{V},\Q)$.
By Theorem \ref{wd}, the two relative theories are determined
by the Gromov-Witten theories of
$V$, $W_1$, and $\widetilde{\proj}(N_1\oplus \oh_{W_1})$ and the classical
restriction maps.

The projective bundle $\proj(N_1|_Z\oplus \oh_Z)$ over $Z$ is a divisor
in $\proj(N_1\oplus \oh_{W_1})$ containing the center $Z\subset D_\infty$ of the blow-up.
The normal bundle of 
$$\proj(N_1|_Z\oplus \oh_Z) \subset \proj(N_1\oplus \oh_{W_1})$$
is the pull-back of $N_2|_Z$.

By repeating the first construction of the proof, we find the
Gromov-Witten theory of $\widetilde{\proj}(N_1\oplus \oh_{W_1})$ is
determined by
the Gromov-Witten theories of
$$\proj(N_1\oplus \oh_{W_1}), \ \proj(N_1|_Z\oplus \oh_Z),$$
and the blow-up of 
$$\proj(N_1|_Z\oplus \oh_Z)\times_Z \proj(N_2|_Z\oplus \oh_Z)$$
along $Z$ embedded as $D_\infty \times_Z D_\infty$.

Finally, the last blown-up variety can be studied via the virtual localization
formula \cite{GP}. The variety
\begin{equation}\label{sx12}
\proj(N_1|_Z\oplus \oh_Z)\times_Z \proj(N_2|_Z\oplus \oh_Z)
\end{equation}
carries a fiberwise $\com^*\times\com^*$-action over $Z$. The action
lifts to the blow-up of \eqref{sx12} along the $\com^*\times\com^*$-fixed
locus $D_\infty \times_Z D_\infty$.
The $\com^*\times \com^*$-action on the blown-up space has
5 fixed loci --- each isomorphic to $Z$. The localization formula
reduces the Gromov-Witten invariants of the blown-up space to
Hodge integrals in the Gromov-Witten theory of $Z$.
\end{proof}

A direct generalization of the proof of Lemma \ref{ss} yields the following
related result. Let $Z\subset V$ be the nonsingular complete intersection of 
$n$ nonsingular divisors
$$W_1, \ldots, W_n\subset V,$$
and let $\widetilde{V}$ be the blow-up of $V$ along $Z$.
The Gromov-Witten theory of $\widetilde{V}$ is determined by the
Gromov-Witten theories of
$$V, \ W_1, \ W_1\cap W_2, \  W_1\cap W_2\cap W_3,\ \ldots,\ \cap_{i=1}^n W_i=Z$$
and the classical restriction maps.

\vspace{+10pt}
\noindent{\bf Conjecture 2.} {\em The Gromov-Witten theory of $\widetilde{V}$
is uniquely and effectively determined by the Gromov-Witten theories of
$V$ and $Z$ and the
restriction map
$H^*(V,\Q) \rarr H^*(Z,\Q).$}
\vspace{+10pt}

\subsection{The Calabi-Yau quintic}
\label{zzwe}
\subsubsection{Notation}
The following notation for curves, surfaces, and 3-folds will be convenient
for the study of the quintic:
\begin{enumerate}
\item[(i)] let $C_{d_1,d_2}\subset \proj^3$ be a nonsingular complete intersection
 of type $(d_1,d_2)$,
\item[(ii)] let $S_d \subset \proj^3$ be a nonsingular surface of degree $d$,
\item[(iii)] let $T_d \subset \proj^4$ be a nonsingular 3-fold of degree $d$.
\end{enumerate}
Finally, let $\proj^3[d_1,d_2]$ be the blow-up of $\proj^3$ along $C_{d_1,d_2}$.

\subsubsection{The quintic scheme}
The calculation scheme for the quintic is diagrammed below by arrows showing the
dependencies of Gromov-Witten theories.
\begin{enumerate}
\item[(a)]  the superscript $*$ denotes simple Gromov-Witten theories,
\item[(b)] the arrow $\stackrel{k}{\longrightarrow}$ denotes the application of
Theorem $k$ in simple or full form,
\item[(c)] the arrow $\stackrel{l}{\longrightarrow}$ denotes the application of
Lemma \ref{ss}.
\end{enumerate}

\noindent{The quintic scheme:}

$$T_5^* \stackrel{3}{\longrightarrow} (T_4,S_4)^*, \ (\proj^3[4,5], S_4),$$
$$(T_4,S_4)^* \stackrel{2}{\longrightarrow} T_4^*, \ S_4,$$
$$T_4^* \stackrel{3}{\longrightarrow} (T_3,S_3)^*, \ (\proj^3[3,4],S_3)$$
$$(T_3,S_3)^* \stackrel{2}{\longrightarrow} T_3^*, \ S_3,$$
$$T_3^* \stackrel{3}{\longrightarrow} (T_2,S_2)^*, \ (\proj^3[2,3],S_2)$$
$$(T_2,S_2)^* \stackrel{2}{\longrightarrow} T_2^*, \ S_2,$$
$$T_2^* \stackrel{3}{\longrightarrow} (T_1,S_1)^*, \ (\proj^3[1,2],S_1)$$
$$(T_1,S_1)^* \stackrel{2}{\longrightarrow} T_1^*, \ S_1,$$
$$ $$
$$(\proj^3[4,5], S_4) \stackrel{2}{\longrightarrow}  \proj^3[4,5], \ S_4,$$
$$\proj^3[4,5] \stackrel{l}{\longrightarrow} \proj^3, \ S_4, \ C_{4,5}, $$
$$(\proj^3[3,4], S_3) \stackrel{2}{\longrightarrow}  \proj^3[3,4], \ S_3,$$
$$\proj^3[3,4] \stackrel{l}{\longrightarrow} \proj^3, \ S_3, \ C_{3,4}, $$
$$(\proj^3[2,3], S_2) \stackrel{2}{\longrightarrow}  \proj^3[2,3], \ S_2,$$
$$\proj^3[2,3] \stackrel{l}{\longrightarrow} \proj^3, \ S_2, \ C_{2,3}, $$
$$(\proj^3[1,2], S_1) \stackrel{2}{\longrightarrow}  \proj^3[1,2], \ S_1,$$
$$\proj^3[1,2] \stackrel{l}{\longrightarrow} \proj^3, \ S_1, \ C_{1,2}. $$
$$ $$

The end points of the quintic scheme are the following absolute Gromov-Witten theories:
$$\proj^3,\ \proj^2,\ S_2, \ S_3, \ S_4, \ C_{1,2}, \ C_{2,3}, \ C_{3,4}, \ C_{4,5},$$ 
all previously determined.

\subsection{The quintic surface $S_5$} \label{exx}
We end the paper with a basic calculation 
for the quintic surface $S_5\subset \proj^3$. Let $K$ denote the canonical class of $S_5$. The
adjunction genus of curves in class $K$ is 6.  The expected dimension of
genus 6 curves of class $K$ is 0.
By Taubes's result equating Gromov and Seiberg-Witten invariants,  
$$\langle 1 \rangle_{6,K}^{S_5}= SW(-P_{S_5}).$$
Since $S_5$ is a minimal surface of general type,
$$SW(-P_{S_5}) = (-1)^{(1+p_{g}(S_5)+q(S_5))} = -1.$$ 
Here, $P_{S_5}$ is the $\mathrm{Spin}^{c}$-structure
induced by the complex structure of $S_5$  and $p_{g}$ and $q$ are 
the geometric genus and irregularity of the surface $S_5$ \cite{morgan}.
We will calculate $\langle 1 \rangle_{6,K}^{S_5}$ directly via our method.

The first step is the degeneration of the quintic to a union of a K3 surface and
a blown-up projective plane along a plane quartic curve:
$$S_5^* \rarr (S_4,C_4)^* , (\proj^2[4,5],C_4).$$
Let $B$ denote $\proj^2[4,5]$, the blow-up of $\proj^2$ in 20 points.
The degeneration formula yields
\begin{equation}\label{quintic}
\langle 1 \rangle^{S_5}_{6,K} = \sum_{g_{i},\beta_{i},\mu} \langle 1|\mu\rangle^{S_4/C_4}_{g_{1},\beta_{1}}
\langle\mu^{\vee}|1\rangle^{B/C_4}_{g_{2},\beta_{2}}
\end{equation}
where the summation is over all possible genus splittings $g_1+g_2+\ell(\mu)-1=6$, class splitting $\beta_1+\beta_2=K$, 
and cohomology-weighted partitions $\mu$.

Let $H,L$ denote the hyperplane classes
on $S_4$ and $B$ respectively and let $E_{i}$ denote the class of the $i^{th}$ exceptional divisor of $B$.

Although all configurations on the right side of \eqref{quintic} can be treated algorithmically, 
most can be easily ruled out without any computation. For example,
configurations with $g_{1}=5,\beta_{1}=H$, require four fixed point condition in $\mu_{1}$. However, since the dimension
of the linear system $|H|$ is 3, the solution set will be empty.  

For configurations with $\beta_{1}=0$, we can show only 
$$\beta_{2} = 5L - \sum_{i=1}^{20} E_{i}$$ is allowed as follows.  
Consider a configuration with 
$$\beta_{2} = 5L - 2E_{1} - E_{3}-\dots - E_{20}.$$
  A curve mapped to $B$ gives a section of $\mathcal{O}_{C}(5)$ on $C$ with divisor $$2p_{1}+p_{3}+\dots+p_{20}.$$  
However, a monodromy argument shows that such a linear equivalence is impossible
for a generic choice of 
$$p_{1}+\dots +p_{20} \in |\mathcal{O}_C(5)|.$$ Finally, for many configurations,
the relative invariant for $(B,C_4)$ is easily seen to vanish.  For instance, 
when 
$g_{1}=3, \beta_2=L$ with $$\mu = \{(2,[p]),(1,1),(1,1)\},$$ the relative invariant for the curve mapped to $B$ is the number 
of tangents to $C_4 \subset \mathbf{P}^{2}$ which pass through
two fixed generic points.

After ruling out these easy cases, there is the single $\beta_1=0$ case discussed above 
and five configurations with $\beta_{1} = H$.
For $g_1= 3$, the possible partitions with $\beta_1=H$ are 
\begin{eqnarray*}
\mu_{1} &= &\{(1,[p]),(1,[p]),(1,1), (1,1)\},\\
\mu_{2} & = &\{(1,\alpha),(1,\alpha^{\vee}),(1,[p]),(1,1)\},\\
\mu_{3} & =& (1,\alpha_{1}),(1,\alpha_{1}^{\vee}),(1,\alpha_{2}),(1,\alpha_{2}^{\vee})\}
\end{eqnarray*}
where $[p] \in H^{2}(C, {\mathbb Z})$ is the Poincare-dual class to a point and $\alpha_{i}, \alpha_{i}^{\vee}$ are
basis elements of $H^{1}(C, {\mathbb Z})$.  For $g_1=4$, the possible partitions are
\begin{eqnarray*}
\nu_{1}&  =& \{(2,[p]),(1,[p]),(1,1)\}, \\
\nu_{2}& = & \{(2,1),(1,[p]),(1,[p])\}.
\end{eqnarray*}
In each of these cases, we have $g_{2}=0$ and $\beta_{2}=L$ for the curve
mapped to $B$.

For each relative invariant of 
$S_4$, we study the associated absolute invariant on $S_4$ via the 
degeneration to
$$S_4 \cup_{C_4} \mathbf{P}(K_{C_4} \oplus \mathcal{O}_{C_4}).$$
We will denote the latter projective bundle by $P$ and the sections
corresponding to the summands $K_{C}$ and $\mathcal{O}_{C}$ by
 $D_{0}$ and $D_{\infty}$.
For $g_1=3$, we obtain the equations
\begin{equation}
\label{sss}
\langle \tau_{0}([p])\tau_{0}([p])\rangle^{S_4}_{3,H}  = \langle 1 | \mu_{1}\rangle^{S_4/C_4}_{3,H} + 
\langle (0)|\tau_{0}([p])\tau_{0}([p])\rangle^{P/D_{\infty}}_{3,[D_{0}]}, 
\end{equation}
\begin{multline}
\label{ssst}
\langle \tau_{0}(\iota_{\ast}\alpha)\tau_{0}(\iota_{\ast}\alpha^{\vee})\tau_{0}([p])\rangle^{S_4}_{3,H}
= \\ \langle 1 | \mu_{2} \rangle^{S_4/C_4}_{3,H} + \langle 1 | \mu_{1} \rangle^{S_4/C_4}_{3,H} + 
\langle (0)|\tau_{0}(\iota_{\ast}\alpha)\tau_{0}(\iota_{\ast}\alpha^{\vee})\tau_{0}([p])\rangle^{P/D_{\infty}}_{3,[D_{0}]},
\end{multline}
and
\begin{multline}
\label{yyy}
\langle \tau_{0}(\iota_{\ast}\alpha_{1})
\tau_{0}(\iota_{\ast}\alpha_{1}^{\vee})\tau_{0}(\iota_{\ast}\alpha_{2})\tau_{0}(\iota_{\ast}\alpha_{2}^{\vee})
\rangle^{S_4}_{3,H} = \\
\langle 1|\mu_{3}\rangle^{S_4/C_4}_{3,H}  + 2 \langle 1 | \mu_{2} \rangle^{S_4/C_4}_{3,H} 
+ \langle 1 | \mu_{1} \rangle^{S_4/C_4}_{3,H} \\
+ 
\langle (0)|\tau_{0}(\iota_{\ast}\alpha_{1})
\tau_{0}(\iota_{\ast}\alpha_{1}^{\vee})\tau_{0}(\iota_{\ast}\alpha_{2})\tau_{0}(\iota_{\ast}\alpha_{2}^{\vee}
)\rangle^{P/D_{\infty}}_{3,[D_{0}]}.
\end{multline}

Since there exist Kahler deformations of the K3 surface $S_4$ with no embedded curves, 
all absolute invariants vanish.  The relative invariants of $(P,C_{\infty})$ are immediately 
calculated via relative localization to yield $1, 0,$ and $0$ respectively.  The genus 3 relative invariants are
$$\langle 1|\mu_{1}\rangle^{S_4/C_4}_{3,H} = -1, \langle 1|\mu_{2}\rangle^{S_4/C_4}_{3,H} = 1, 
\langle 1 |\mu_{3}\rangle^{S_4/C_4}_{3,H}=-1.$$
The same procedure for the two genus 4 invariants shows they both vanish.

It remains to compute the relative invariants of the pair $(B,C_4)$.  
For the relative invariants $\langle \mu_{i}^{\vee}|1\rangle^{B/C_4}_{0,L}$,
we argue directly.  The moduli space $$\overline{M}_{0,4}(B/C_4, (1,1,1,1))$$ is a blowup of the two 
dimensional space
$$Z = \{(p_{1},p_{2},p_{3},p_{4}) \in C \times C\times C\times C \ |\  \mathcal{O}_{C}(\sum p_{i}) \cong K_{C}\}.$$  
Since the relative 
conditions are pulled back from $C^{4}$, we only need to cap these classes with the fundamental class of $Z$ in 
$H_{2}(C^{4},{\mathbb Z})$.  The latter class is given by a 
degeneracy locus calculation, $[Z]$ equals the second Chern class of 
the rank 4 bundle whose fiber at $(p_{1},p_{2},p_{3},p_{4})$ is $\bigoplus K |_{\sum p_{i}}$.  The invariants are
$$\langle \mu_{1}^{\vee}|1\rangle^{B/C_4}_{0,L}= 1, \langle \mu_{2}^{\vee}|1\rangle^{B/C_4}_{0,L} = 1, 
\langle \mu_{3}^{\vee}|1\rangle^{B/C_4}_{0,L} = 1.$$

The remaining contribution is the relative invariant for $(B,C_4)$  
where the entire genus 6 curve is mapped to $B$ with homology class $\beta = 5H - \sum_{i=1}^{20} E_{i}$.
The corresponding absolute invariant is
\begin{equation}\label{ccffw}
\langle 1 \rangle^{B}_{6,\beta} = 1.
\end{equation}
The invariant \eqref{ccffw}
 is computed by observing that, given two generic curves $C_{4}$ and $C_{5}$ of degree 4 and 5 in $\mathbf{P}^{2}$,
 $C_{5}$ is the unique quintic curve through the 20 points of $C_{4}\cap C_{5}$.  
Degenerating along $C_4$ gives the following relation
between the relative invariants for $(B,C_4)$ and $(P,D_{0})$.
\begin{equation}\label{blowup}
1 =\langle 1 \rangle^{B}_{ 6,\beta} = \sum_{g_{i},\beta_{i},\mu} \langle 1|\mu\rangle^{P/D_{0}}_{g_{1},
\beta_{1}}\langle\mu^{\vee}|1\rangle^{B/C_4}_{g_{2},\beta_{2}},
\end{equation}
where again we sum over all configurations.  

The situation is identical to our degeneration of $S_5$ with $P$ in place of $S_4$.  
In particular, the same arguments allow us to reduce to the same five configurations with
$g_{1}=3$ or $4$, in addition to the configuration where the entire curve is mapped to $B$. 
 Moreover, the relative invariants are computed by the same set of equations
(\ref{sss}-\ref{yyy}) with $P$ instead of $S_4$. 
The only difference is the absolute invariants on $P$ no longer vanish identically and must be computed directly by localization.  
Omitting the details, we find for genus 3, the relative invariants are
$$\langle 1|\mu_{1}\rangle^{P/D_{0}}_{3,[D_0]} = 7, \
\langle 1|\mu_{2}\rangle^{P/D_{0}}_{3,[D_0]} = -3, \
\langle 1 |\mu_{3}\rangle^{P/D_{0}}_{3,[D_0]}=1.$$
For genus 4, the relative invariants again vanish.

Finally, if we subtract equation \eqref{blowup} from equation \eqref{quintic}, we have 
\begin{eqnarray*}
\langle 1\rangle^{S_5}_{6,K} & = &1+ \sum_{i} \left(\langle 1|\mu_{i}\rangle^{S_4/C_4}_{3,H} - 
\langle 1|\mu_{i}\rangle^{P/D_{0}}_{3,[D_{0}]}\right)
\langle \mu_{i}^{\vee}|1\rangle^{B/C_4}_{0,L}
\\
& = & 1+(-1-7)\cdot 1 + 3(1-(-3))\cdot 1 + 3(-1-1)\cdot 1\\ & = &-1
\end{eqnarray*}
where the factors of 3 arise from the different combinations of odd cohomology conditions $\alpha \in H^{1}(C,{\mathbb Z})$.


\vspace{+10 pt}
\noindent
Department of Mathematics \\
Princeton University \\
Princeton, NJ 08544, USA\\
dmaulik@math.princeton.edu \\

\vspace{+10 pt}
\noindent
Department of Mathematics\\
Princeton University\\
Princeton, NJ 08544, USA\\
rahulp@math.princeton.edu

\end{document}